\documentclass[12pt, a4paper]{article}
\usepackage{latexsym} 
\usepackage[all]{xy}
\usepackage{amsmath}
\usepackage{amsfonts}
\usepackage{amssymb}
\usepackage{mathrsfs}

\newtheorem{thm}{Theorem}[section]
\newtheorem{prop}[thm]{Proposition}

\newtheorem{df}[thm]{Definition}
\newtheorem{cor}[thm]{Corollary}
\newtheorem{conj}[thm]{Conjecture}

\def\k{{\kappa}}
\def\x{\times}

\def\shO{{\cal O}}
\def\shG{{\cal G}}
\def\scG{{\mathscr G}}
\def\shU{{\cal U}}
\def\shA{{\cal A}}
\def\cG{{\cal G}}

\def\bC{{\mathbb C}}

\def\bR{{\mathbb R}}

\def\bV{{\mathbb V}}

\def\bT{{\mathbf T}}

\newcommand{\func}[1]{{\stackrel{#1}{\longrightarrow}}}

\newcommand{\pf}{{\bf Proof. }}
\newcommand{\qed}{{$\square$}}
\newcommand{\tn}[1]{\textnormal{#1}}
\newcommand{\cat}[1]{(\tn{#1})}
\newcommand{\remark}{\noindent{\bf Remark. }}

\newcommand{\st}[1]{{\mathfrak #1}}
\newcommand{\hst}[1]{\hat{\mathfrak #1}}

\newcommand{\noi}{\noindent}

\title{Nonabelian Mixed Hodge Structure on Brill-Noether Stacks}
\author{Vilislav Boutchaktchiev}
\date{October 1, 2013}

\begin{document}

\maketitle
\begin{abstract}
A Brill-Noether stack is an algebraic very presentable stack whose
homotopy type has two nontrivial homotopy groups. We consider one
with a fundamental group  --- a reductive algebraic group-scheme $S$
and one higher homotopy group, represented by a vector space $V$. The
homotopy type also defines an action of $S$ on $V$. This stack is used
as coefficient space for nonabelian cohomological space on a smooth
algebraic variety $X$.  

We define nonabelian MHS on cohomological spaces of this type in the
context of the work of C. Simpson related to MHS on the space of local
systems. It is defined via an action of the multiplicative complex
group on a appropriately chosen category. The exhibited structure in
fact generalizes Simpson's work. Allowing a more general type of
coefficient stack. 

Furthermore, the so-defined MHS, when considered at a vicinity of an
object, which remains fixed under the structural action of $bC^*$,
produces the local MHS on Brill-Noether Stacks as defined in our
earlier work.   

The nonabelian mixed Hodge structure on a Brill-Noether stack is a
example of the mixed Hodge structure on a schematic homotopy type,
studied by Katzarkov, Pantev and Toen. It has the advantage, due to
the relative simplicity of the coefficient stack, that it could be
locally written out in terms of iterated integrals. 
\end{abstract}

\tableofcontents

\section{Introduction}
Let $X$ be a smooth complex algebraic curve. We study the geometry of $X$
through introducing a Hodge structure on the nonabelian cohomology
space $Hom(X,T)$, where the coefficient space, $T$,  is an algebraic
stack. To simplify the problem, one can restrict $T$ to be a very
presentable geometric $n$-stack, which, in the category of topological
spaces,
is analogous to considering a homomorphism with to an $n$-truncated 
CW-complex. The case of $X$ being a curve of positive 
genus is additionally simplified by the fact that $X$ is a
topologically $K(\pi_1,1)$-space.

This problem was studied before by Simpson, in e.g. \cite{Si:1, Si:2},
where he defined a Hodge 
structure on the stack $Hom(X,\k(G,1))$, which parametrizes $G$-local
systems on $X$ up to homotopy equivalence. Simpson proves that each
$\bC$-valued point corresponds to a Higgs bundle with a holomorphic connection
$(E,\theta)$ and, on the other hand, to a representation
$\sigma:\pi_1(X,x)\to G$. The ring of functions $\shO(E)$ reflects the
local geometry of $Hom(X,k(G,1))$ at the point $(E,\theta)$ and
Simpson proves, that his Hodge structure restricts to a mixed Hodge
structure on $\shO(E)$, which, under some additional restrictions,
coincides with the MHS defined by Hain \cite{Ha:1}.

In this paper we suggest a MHS on $Hom(X,T)$, where $T$ is a {\it
  Brill-Noether stack}, i.e., a stack which has only 2 nontrivial
homotopy groups: $\pi_1=G$ and  $\pi_n=V$ and $\pi_1$ acts on
$\pi_n$ by  a chosen   representation $\rho:G\to GL(V)$. For that,
to a $\bC$-valued point of $Hom(X,T)$, which parametrizes a triple
  $(E,\theta,\eta)$ of a flat $G$-bundle, holomorphic connection on
  $E$ end a cohomology class $\eta\in H^n(X,E\x^GV)$,
we define a MHS on a proalgebraic group $\shG_\rho$
and  
we define an action of $\shG$ on $\shG_\rho$, which reflects the
Whitehead product on the homotopy type of $T$. The resultant MHS on
$Hom(X,T)$ is comprise by 2 ingredients --- Simpson's nonabelian
MHS on $Hom(X,\k(G,1))$ an the usual abelian Hodge theory on the
cohomology group $H^n(X,E\x^GV)$. 

In section~\ref{sec:1} we provide some background information about the
theory of stacks. 

In section~\ref{sec:3} contains the main result,
Theorem~\ref{thm:hodge}. In a previous work we have defined local
mixed Hodge structure on a Brill--Noether stack. Here we justyfy the
name {\it local}, proving that at a given complex point of $Hom(X,
BN)$, fixed under the defining $\bC^*$-action for the Hodge structure,
the the resultant filtration on the local algebra is the same as
defined before. 

In section~\ref{sec:3} we discuss, also, some questions of
compatibility with existant 
literature give a few conjectures which    demonstrate that
our 
approach can be generalized in durther study to define 
MHS  on $Hom(X,T)$, for any algebraic smooth variety $X$ and any
topological sheaf $T$.

\section{The Brill-Noether stack}\label{sec:1}
\subsection{The Eilenberg-MacLane stacks $\k(G,n)$.}
In the this section we recall the definitions and basic properties of 
$\k(G,n)$, following \cite{Si:1} and \cite{To:1}.  

We denote by $\bT$ the category of topological spaces. For us a
{\it $n$-stacks} on the site $\bT$ will be $n$-truncated   
objects in the simplically enriched category $St(\bT)=LSPr(\bT)$,
which is obtained in \cite{To:1} by taking the category of simplicial 
presheaves on $\bT$, applying simplicial localization with respect to
the local equivalences and introducing simplicial structure on the
arrows. According to \cite{Si:1}, these stacks correspond to {\it very
  presentable geometric $n$-stacks} on $\cat{Aff}_{\textit{\'et}}$,
via a theorem of GAGA type. This means, that  $\pi_1(T,t)$ is
represented by an algebraic group-scheme of finite type, and
$\pi_j(T,t)$ is represented by a vector space for $j>1$.

If $G$ is a group-sheaf, represented by an algebraic scheme $G$,  a
simple example of a very presentable stack, one can define $\k(G,n)$, so
that $\pi_j(\kappa(G,n))$ is trivial except $\pi_j=G$. ($G$ must be
abelian if $n>1$.)  

\remark For $n=1$, this is the classifying stack $BG$ of $G$-torsors. 
For $n>1$, $\k(G,n)$ can be constructed as a topological realization
of certain sheaf  of Lie algebras.

The stacks $\k(G,n)$ are classifying for the group $G$ in the sense of
the following proposition:
\begin{prop}[Simpson, \cite{Si:2}]\label{prop:1}
For each $X\in\bT$, $Hom(X, \k(G,n))$ is a stack, such that
\[\pi_j(Hom(X, \k(G,n))=H^{n-j}(X,G).
\]
\end{prop}
This proposition gives a common expression for two separate notions:
\begin{enumerate}
  \item When $n=1$ and $G$ is any sheaf of groups, $H^1(X,G)$ is the stack
    parametrizing the $G$-torsors on $X$, i.e., the pairs 
    $(E,\theta)$, where $E$ is a     flat $G$-bundles over $X$, and
    $\theta$ is a flat holomorphic connection.
  \item If  $n\ge0$, and $G=V$ is a $\bC$-vector space, $H^i(X,V)$ is
    the usual abelian cohomology and the $\bC$-valued points of 
$Hom(X, \k(V,n))$ correspond to pairs
$(E,\eta)$, where $E$ is isomorphic to the trivial vector bundle
$\bV=V\times X\to X$ and $\eta\in H^n(X,\bV)$. In this way, we can
introduce a structure of $n$-stack on $H^n(X,\bV)$.
\end{enumerate}

As an extension,  $Hom(X,\k(V,1))$ can be viewed as both $BV(X)$ or
$H^1(X,\bV)$.

%

\begin{prop}
There is an isomorphism:
\[BV(X)\cong H^1(X,\bV)
\]
\end{prop}
\pf
For each flat $V$-bundle the flat connection $\theta$ can be considered
as a cohomological class from $H_{DR}^1(X,\bV)$, which is an
expression of the usual, abelian Hodge theory in degree 1.\qed

\subsection{Postnikov towers}

Very presentable geometric stack can be built as Postnikov towers of
$\k(G,n)$-s. The simplest example is the {\it Brill-Noether stack},
denoted $\k(G,\rho,n)$, whose $\pi_1=G$, $\pi_n=V$, the
representation  $\rho:G\to GL(V)$ gives the action of $\pi_1$ on
$\pi_n$ and the rest of the homotopy groups being
trivial. $\k(G,\rho,n)$ is a fibration over $\k(G,1)$:
\begin{equation}\label{eq:1}
\xymatrix 
{
\k(V,n)\ar[r]^i&\k(G,\rho,n)\ar[d]^\tau\\
&\k(G,1)
}
\end{equation}
The map $\tau$ is truncation and $i$ is inclusion.

For each $X\in\bT$, $Hom(X,\k(G,\rho,n))$ parametrizes triples 
$(E,\theta, \eta)$ of a $G$-bundle $E$, flat connection $\theta$ on
$E$ and a cohomology class $\eta\in H^n(X,E\times^G V)$ (cf. \cite{Si:1}).
The diagram (\ref{eq:1}) induces a smooth fibration of geometric $n$-stacks
\[Hom(X,\k(G,\rho,n))\func{\phi} Hom(X,\k(G,1))
\]
Over each $\bC$-point $(E,\theta)$ of $Hom(X, \k(G,1))$, the fiber is
$H^n(X,E\times^G V)$.

\noi{\bf Examples:}

\noi{\bf A.} For any algebraic group $H$ one can introduce as a special case of
Brill-Noether stack ${\cal G}_H$ as an extension:
\begin{equation}\label{eq:2}
\xymatrix 
{
\k(Z(H),2)\ar[r]^i&{\cal G}_H\ar[d]^\tau\\
&\k(Out(H),1)
}
\end{equation}
$Hom(X,\cG_H)$ is then a 2-stack, parametrizing the 1-stacks $\st{X}\to
X$ which are locally equivalent to $Hom(X,\k(G,1))$.

\noi{\bf B.} Let G be an complex algebraic group, which splits into 
a semidirect product $G=V\rtimes S$ of a reductive
subgruip $S$ and an abelian group $V$ by the map $\rho S\to Aut(V)$.
\begin{prop}\label{prop:1.2.1}
  There is an equivalence of categories
\[Hom(X,\k(G,1))\cong Hom(X,\k(S,\rho,1))
\]
\end{prop}
\qed
\subsection{Universal Fibrations}
In this paragraph we make a clarification on the terminology.

If we consider the $\k(V,n)$ as an appropriate quotient stack we can
restate Proposition \ref{prop:1} in the following form: 
\begin{prop}
  There is an universal bundle $*\to\k(V,n)$ from which every
  $n$-torsor with structural group $V$ $E\to X$ is a pullback via
  certain morphism $X\to\k(V,n)$.
\end{prop}
\begin{cor}
  There is an universal fibration $*\to\k(G,\rho,n)$, in the
  category of triples $(E,\theta,\eta)$, as above. 
The Higgs bundles which have $h^n(X,E\x^GV)\ge 1$ form an open locus in
  $Hom(X,\k(G,\rho,n))$ consisting from the non-trivial fibers in
  (\ref{eq:1}).
\end{cor}
\pf The suggested universal fibration is a pull-back from the
universal $G$-bundle $*\to\k(G,1)$ through the canonical morphism
(\ref{eq:1}). Every map $f:X\to\k(G,\rho,n)$ corresponds to a Higgs bundle
$(E,\theta)$, which is a pull back of the universal $G$- bundle by the
composition of $f$ with (\ref{eq:1}). By universality of pull-back,
the  corresponding to $f$ triple, is a pull-back from $E$. \qed
\bigskip

\remark From the usual abelian Hodge theory follows that:
$$H^n(X,E\x^GV)\cong H^0(X,E\x^GV\otimes\Omega^n)= \overline{H^0(X,E\x^GV})
$$ 

\section{Nonabelian mixed Hodge structure}\label{sec:3}
\subsection{MHS on $Hom(X,\k(G,1))$}\label{ssec:3.1}
In this section we briefly remind some definitions and facts of  
the nonabelian Hodge theory of $Hom(X,\k(G,1))$, as developed in
\cite{Si:3}.

The classical (abelian, pure) Hodge structure can be defined using an
$C^*$ action on the underlying vector space $A_\bC$, which is usually
a cohomological algebra, and deriving the filtrations from the
eigenspaces of the action.  
More formally (cf. \cite{Si:3}),
\begin{df}
  A non-abelian Hodge structure is a triple $(\scG,\Gamma, \mu)$,
  consisting of
  \begin{enumerate}
  \item  an affine group scheme (pro-algebraic group) $\scG$, defined
    over $\bR$ 
  \item a finitely generated subgroup $\Gamma\subset\scG_\bR$, which is
    Zariski dense in $\scG$.
  \item a group action $\mu:U(1)\x\scG_\bR\to \scG_\bR$ via
    homomorphism of the proalgebraic group, such that, the map
    $U(1)\x\Gamma\to \scG^{an}$ is continuous and the element $C=-1$
    in $U(1)$ is a Cartan involution of $\scG_\bR$.
  \end{enumerate}
\end{df}

We consider the stack $Hom(X,\k(G,1))$ as a cohomological space based
on \ref{prop:1}. 

The nonabelian Hodge structure is defined by introducing an action on
the objects of the 
category $Hom(X,\k(G,1))(U)$. Each object in there  parametrizes a
pair $(E,\theta)$, where $E$ is a Higgs bundle on $X\x U$ and $\theta$
is a Higgs field on $E$, i.e., a holomorphic 1-form with coefficient
in $E$. Simpson defines a $C^*$ action on $Hom(X,\k(G,1))(U)$ by
\[
t:(E,\theta)\mapsto(E,t\theta). 
\]

In this case the group $\scG$ is the pro-algebraic completion of the
fundamental group $\pi_1(X,x_0)$ and the action $U(1)$ is induced from
the described geometric action of  $\bC^*$.

To compare this theory with the MHS found by Hain, one  considers 
any  $\bC$-valued point of $Hom(X,\k(G,1))$ as a representation
$\sigma:\pi_1(X,x_0)\to G$,  which gives a relative prounipotent
completion $\shG$. 

Recall, that $\sigma$ is called to be {\it of Hodge type} if ti comes
from an admissible complex variation of Hodge structure. 
Simpson proves the following proposition:
\begin{prop}
  $\sigma$ comes from a complex variation of Hodge structure if, and
  only if, it is a fixed point of the $\bC^*$ action on
  $Hom(X,\k(G,1))$.
\end{prop}

In this case, 
there is an induced action on the Lie algebra of the progroup
$\shG$ and hence, 
 a decreasing filtration $F^\bullet$ on ${\mathrm{Lie}\,}\shG$, defined by:
\[F^p=\oplus_{r\ge p} H^p, 
\]
where $H^r\subset {\mathrm{Lie}\,}\shG$ is the subspace  of elements
$h$, such that $t(h)=t^r$ for all $t\in\bC$. The weight filtration is,
again as in (\ref{pr:2.1.3}), defined  in the terms of the lower
central series of the Lie algebra of the unipotent radical of $\shG$.
\begin{thm}[Simpson \cite{Si:3}]
  For the trivial representation $\pi_1(X,x_0)\to \{1_G\}$ (which is
  definitely fixed) the mixed
  Hodge structure on the prounipotent completion $\shG$, defined
  above, coincides with Hain's (as described in (\ref{pr:2.1.3})).
\end{thm}
\subsection{Proposed Hodge structure on the Brill-Noether stack}
Let $X$ be a smooth algebraic curve and let $E$ be any quasi-coherent
sheaf of rings on $X$. The d.g.a 
$\Omega^\bullet(X;E)$ has a natural
$\bC^*$ action: 
\begin{equation}\label{eq:action}
t:\omega\mapsto t^r\omega,\qquad \tn{if } \omega\in\shA^{r,s}(X;E)
\end{equation}
This action induces one on the cohomology
d.g.a. $H^\bullet(X;\Omega^r(X;E))$.

Further, for each locally free $E$-module $M$, there is a
$\bC^*$-action on $\shA^\bullet(X;M)$ and hence on  
$H^\bullet(X;\Omega^r(X;M))$.

Suppose, that $T=\k(G,\rho,n)$ is a Brill-Noether stack, such that
$\rho$ is a representation of Hodge type\footnote{The assumption that 
  $\rho:G\to GL(V)$ must be a homomorphism of Hodge type is not
  essential. It can be 
  omitted by considering mixed twistor structures. Nonetheless, even if
 in this case, some special points of $Hom(X,\k(G,\rho,n))$ will carry
 Mixed Hodge structure, due to additional symmetries that they may
 have. This approach shall be discussed in the future.}.
The
$\bC$-points of $Hom(X,T)$ are expressed in the terms of
(eventually, nonabelian) cohomology classes of 
 $H^\bullet(X;\Omega^r(X;E))$ and $H^\bullet(X;\Omega^r(X;M))$. This
 defines an $\bC^*$-action on  $Hom(X,T)$. 
\begin{conj}
There is a non-abelian MHS on $Hom(X,\k(G,\rho,n))$, in the sense of
\cite{K-P-T}, defined by the action \eqref{eq:action}.

\end{conj}
Let $[b]$ be a $\bC$-point of $Hom(X,T)$, which is fixed under
that action.  Then there are  induced actions on the relative 
prounipotent completion $\shG$ of the $\pi_1(X,x_0)$, correspondent to
$E$ and on all the relevant quotients $\shG_\rho$ and therefore on the
correspondent Lie algebras $\hst{g}$ and $\hst{g}_\rho$.

\begin{thm}\label{thm:hodge}
  The MHS's on $\hst{g}$ and on $\hst{g}_\rho$, with Hodge
  filtration induced from the $\bC^*$ action and weight filtration
  coming from the unipotent radical of $\shG$ coincide with those
  described in section \ref{sec:2}. The natural action
$\hst{g}\otimes\hst{g}_\rho\to\hst{g}_\rho$, induced
  by the Whitehead product of the homotopy type of $T$.
is  is given by Theorem \ref{thm:3.2.1}.
\end{thm}
\qed
\subsection{Further conjectures}
The described mixed Hodge structure on $Hom(X,T)$, where $X$ is a
$K(\pi_1,1)$ complex space and $T$ is a Brill-Noether stack, suggests
the following generalizations:

\noi{\bf A:} Let $T$ to be a perfect complex, i.e., the
homotopy type of $T$ is analogous to a simply-connected topological
space. If $X$ is still a topological $K(\pi_1,1)$ space this means
that the nonabelian part of the MHS on $Hom(X,T)$ will be
trivial. Locally, we expect the following:
\begin{conj}\label{conj:perfect}
   Suppose $T$ is constructed as a Postnikov tower 
\[T=T_n\to \dots\to T_i\to T_{i-1}\to\dots\to F_2\to T_0=*
\]
and on each step we have a fibration:
\[
\xymatrix 
{
\k(V_i,i)\ar@{^{(}->}[r]&T_i\ar[d]^{\tau_i}\\
&T_{i-1}
}
\]
Then for each $i$ there exists a proalgebraic group $\shU_i$ (which
will be unipotent), with a mixed Hodge structure, corresponding to the
cohomologies of the of the trivial principal bundle. The $\shU_i$-s form
a mixed Hodge complex with  maps
$\st{u}_i\otimes\st{u}_j\to\st{u}_{i+j-1}$ coming from Whitehead
product on $T$.
\end{conj}
\noi{\bf B:} Let $T$ to be a very presentable $n$-stack, which is a
geometrical realization of a general topological presheaf with $\pi_1=G$
and $\pi_i=V_i$ for $i=2\dots n$. Then for each $i=2...n$ there
is a smooth fibration of algebraic stacks: $T\to\k(G,\rho_i,i)=F_i$, where
$\rho_i:G\to GL(V_i)$ is a part of the definition of $T$.
\begin{conj}
  There exist proalgebraic groups $\shG_i=\shU_i\rtimes G$, (for
  $i=2\dots n$) with MHS on each of them, coming from the corresponding
  Brill-Noether stack $F_i$. They form a MH complex with structure
  maps from (\ref{conj:perfect}).
\end{conj}
\noi{\bf C:} Let $X$ to be any $\bC$-space. Then we can represent it
as a fibration
\[
\xymatrix 
{
X_1\ar@{^{(}->}[r]&X\ar[d]\\
&K(\pi_1,1)
}
\]
where $X_1$ is simply-connected. 
Then MHS on $Hom(X,T)$ should combine the MHS on $Hom(K(\pi_1,1),T)$
from before and $Hom(X_1,T)$ which can be defined by modifying
the MHS of Morgan and Hain on the cohomology algebra of
simply-connected manifold. Locally, on the MH complex language, this means that
we will have a HS for each $Hom(\pi_i(X),\pi_j(T))$, some abelian and
some not, related between each other through Whitehead product in $T$
or $X$.

A reasonable definition of a MHS in all of the cases {\bf A}, {\bf B},
and {\bf C} should respect the MHS on the schematization of  the
homotopy type,  according to \cite{K-P-T}.

\end{document}